\documentclass[11pt]{amsart}
\usepackage{amssymb}
\input xypic
\xyoption{all}
\usepackage[
            pdftex,
           colorlinks=true,
            urlcolor=blue,       
            filecolor=blue,     
            linkcolor=blue,       
            citecolor=blue,         
            pdftitle= {isogenies between Abelian surfaces},
            pdfauthor={Igor Dolgachev, David Lehavi},
            pdfsubject={Abelian surfaces},
            pdfkeywords={}
            pagebackref,
            pdfpagemode=None,
            bookmarksopen=true]
            {hyperref}
\usepackage{hyperref}
\usepackage{color}

\newcommand{\BZ}{{\mathbb{Z}}}
\newcommand{\BK}{{\mathbb{K}}}

\newcommand{\BC}{{\mathbb{C}}}
\newcommand{\BF}{{\mathbb{F}}}

\newcommand{\BP}{{\mathbb{P}}}

\newcommand{\BG}{{\mathbb{G}}}
\newcommand{\BV}{{\mathbb{V}}}
\newcommand{\BT}{{\mathbb{T}}}

\newcommand{\St}{{\textup{St}}}

\newcommand{\gT}{\Theta}
\newcommand{\calM}{{\mathcal{M}}}
\newcommand{\calA}{{\mathcal{A}}}
\newcommand{\calO}{{\mathcal{O}}}

\newcommand{\calK}{{\mathcal{K}}}

\newcommand{\calL}{{\mathcal{L}}}
\newcommand{\calH}{{\mathcal{H}}}
\newcommand{\calX}{{\mathcal{X}}}
\newcommand{\calB}{{\mathcal{B}}}
\newcommand{\calC}{{\mathcal{C}}}
\newcommand{\calS}{{\mathcal{S}}}

\newcommand{\calZ}{{\mathcal{Z}}}

\newcommand{\Pic}{\mathrm{Pic}}
\newcommand{\Hom}{\mathrm{Hom}}
\newcommand{\jac}{\mathrm{Jac}}

\newcommand{\Kum}{\mathrm{Km}}
\newcommand{\NS}{\mathrm{NS}}
\newcommand{\cha}{\mathrm{char}}
\newcommand{\half}{\frac{1}{2}}
\newcommand{\sm}{\smallsetminus}
\newcommand{\Jac}{\mathrm{Jac}}
\newcommand{\la}{\langle}
\newcommand{\ra}{\rangle}

\newcommand{\PSp}{\textup{PSp}}
\newcommand{\Sp}{\textup{Sp}}
\newcommand{\Th}{\textup{Th}}
\newcommand{\Hess}{\textup{Hess}}
\newcommand{\Or}{\textup{O}}

\newcommand{\beq}{\begin{equation}}
\newcommand{\eeq}{\end{equation}}
\newcommand\thetafun[2]{\vartheta\left[\begin{smallmatrix}#1\\ #2\end{smallmatrix}\right]}
\theoremstyle{plain}
\newtheorem{lma}{Lemma}[section]
\newtheorem{thm}[lma]{Theorem}
\newtheorem{prp}[lma]{Proposition}
\newtheorem{cor}[lma]{Corollary}
\theoremstyle{definition}

\newtheorem{rmr}[lma]{Remark}

\begin{document}


\title[Isogenous abelian surfaces]{On isogenous principally polarized abelian surfaces}
\author{I. Dolgachev}
\address{Department of Mathematics, University of Michigan,
 530 Church St, Ann Arbor 48109 USA}
\email{idolga@umich.edu}
\author{D. Lehavi}
\address{Correlix Ltd. 6 Galgaley Haplada St,
Herzelia Pituah, Israel 46733}
\email{dlehavi@gmail.com}
\date{}
\dedicatory{To Roy Smith}
\begin{abstract}
We study a relationship between two genus $2$ curves whose jacobians are isogenous with kernel equal to a maximal isotropic subspace of $p$-torsion points with respect to the Weil pairing. For $p = 3$ we find an explicit relationship between the set of Weierstrass points of the two curves extending the classical results of F. Richelot (1837) and G. Humbert (1901) in the case $p = 2$.  
\end{abstract}
\maketitle

%
\section{Introduction}\label{Sintro}
%

Let $A = \BC^g/\BZ^g+\tau\BZ^g$ be a complex algebraic torus with the period matrix $\tau$ in the Siegel space $Z_g$. Replacing $\tau$ by $p\tau$ for some integer $p$ defines an isogenous torus with the kernel of isogeny equal to $F = (\BZ/p\BZ)^g$. In the case of elliptic curves this construction leads to the theory of modular equations and, for $p = 2$, to the  Gauss algebraic-geometrical mean. The  construction can be made independent of the choice of $\tau$ and, in fact, can be defined for any ordinary abelian principally polarized variety $A$ over an algebraically closed field $\BK$ of arbitrary 
characteristic. We take a maximal isotropic subspace $F$ in the group $A[p]$ of $p$-torsion points with respect to the Weil pairing defined by the principal polarization. Then the quotient $B = A/F$ is a principally polarized abelian variety. When $\BK = \BC$, the isomorphism class of $B$ can be defined by the period matrix $p\tau$, where $\tau$ is a period matrix of $A$.

In the case where $A = \Jac(C)$ is the jacobian variety of a smooth algebraic curve of genus $g$ over $\BK$, one may ask whether $B = \Jac(C')$ for some other curve $C'$ of genus $g$, and if so, what is the precise relationship between the moduli of $C$ and $C'$. For $g = 1$ the answer is given by exhibiting an explicit modular equation relating the absolute invariants of the two elliptic curves (see a survey of classical and modern results in \cite{BB}). In the case $p = 2$ this was done by Gauss in his work about the algebraic-geometrical mean (see \cite{Cox}). For  $g = 2$ and $p = 2$, the explicit geometric moduli relationship between the two curves of genus $2$ was found by Richelot \cite{Ric}  and Humbert \cite{Hum} (see a modern account in \cite{BM}).  It was extended  to the case $g= 3$  
by Donagi-Livn\'e \cite{DoLi} and Lehavi-Ritzenthaler \cite{LR}.  

 In the present paper we study the case $g = 2$ and $p > 2$ and assume that the ground field $\BK$ is an algebraically closed field of characteristic $\ne 2$.  Our main result is the following

\begin{thm}\label{Tmain}
Let $C$ be a smooth projective genus $2$ curve. Let $p$ be an integer coprime to the characteristic of $\BK$, and $F$ be a maximal isotropic subgroup of $\Jac(C)[p]$ with respect to the Weil pairing defined by the natural principal polarization of the jacobian.  For any $e,-e\in F\sm \{0\}$ let $\{x_e,y_e\}$ be a unique, up to the hyperelliptic involution $\iota$, pair of points on $C$ such that $x_e-y_e = \pm e$. Let $\phi:C\to R_{2p} \subset \BP^{2p}$ be the degree two map onto a rational norm curve given by the linear system $|2pK_C|^\iota$. Let $(c_e,d_e)$ be the images of the pairs $(x_e,y_e) $ in $R_{2p}$ and $\ell_e = \overline{c_e,d_e}$ be the corresponding secant lines of $R_{2p}$. There exists a unique hyperplane $\calH$ in $\BP^{2p}$ containing the images $w_1,\ldots,w_6$ of the six Weierstrass points  such that the intersection points of $\calH$ with the secants $\ell_e$ are contained in  a subspace $L$ of $\calH$ of codimension $3$. The images of the points $w_i$ under a projection from $L$ to $\BP^3$ are contained in a conic (maybe reducible), and the double cover of the conic ramified at these points is a stable curve $C'$ of arithmetic genus $2$  such that $\Jac(C') \cong \Jac(C)/F$.
\end{thm}

In the case $\BK = \BC$ and $p = 3$ we give an effective algorithm for determining the curve  $C'$ in terms of $C$. The assertion of the theorem is also valid in the case $p =\cha~ \BK$ if we assume that 
$F = \Jac(C)[p]\cong (\BZ/p\BZ)^2$, i.e. $\Jac(C)$ is an ordinary abelian variety. In this case the algorithm for finding $C'$ can be made even more effective. 

Note that G. Humbert also considered the case $g = 2, p = 3$ in \cite{Hum} but his solution is different from ours and cannot be made effective (see Remark \ref{humbert}).
%
\section{Preliminaries}
%
\subsection{Polarized abelian varieties}\label{weil}  
Let $A$ be a $g$-dimensional  abelian variety over an algebraically closed  field $\BK$. Let $\calL$ be an  invertible sheaf on $A$ and $\pi:\BV(\calL)\to A$ be the corresponding line bundle, the total space of $\calL$. One defines the \emph{theta group scheme} $G(\calL)$ whose $S$-points are lifts of translation automorphisms $t_a, a\in A(S),$ of $A_S = A\times_\BK S$ to automorphisms of $\BV(\calL)_S$. The theta group scheme fits in the canonical central extension of group schemes 
\[
1\to \BG_m \to G(\calL) \to K(\calL) \to 1,
\]
where $K(\calL)(S)$ is the  subgroup  of $A(S)$ of translations which send $\calL_S$ to an isomorphic invertible sheaf on $A_S$. The extension is determined by the \emph{Weil pairing}
\[e^{\calL}:K(\calL)\times K(\calL)\to \BG_m\]
defined by the commutator in $G(\calL)$. A subgroup $K$ of $K(\calL)$ is isotropic with respect to the Weil pairing if and only if the extension splits over $K$. 

From now we assume that $\calL$ is ample. In this case $K(\calL)$ is a finite group scheme and the Weil pairing is non-degenerate. Recall that the algebraic equivalence class of an ample invertible sheaf on  $A$ is called a \emph{polarization} on $A$. An abelian variety equipped with a polarization is called a \emph{polarized abelian variety}. 

Let $A^\vee$ be the dual abelian variety representing the connected component of the Picard scheme of $A$. Any invertible sheaf  $\calL$ defines a homomorphism of abelian varieties 
\[
\phi_\calL:A\to A^\vee, a\mapsto t_a^*(\calL)\otimes \calL^{-1}.
\]
The homomorphism depends only on the algebraic equivalence class of $\calL$ and its kernel is isomorphic to the group $K(\calL)$.  In particular, $\lambda$ is an isogeny if and only $\calL$ is ample. We say that $\calL$ defines  a {\em principal polarization}  if $\phi_\calL$ is an isomorphism. This is also equivalent to $\calL$ being ample and $h^0(\calL) = 1$. 

The proof of the following proposition can be found in \cite{Mu2}, \S 23.

\begin{prp} Let $\lambda:A\to B$ be a separable isogeny of abelian varieties. There is a natural bijective correspondence between the following sets
\begin{itemize}
\item the set of isomorphism classes of invertible ample sheaves $\calM$ such that $\lambda^*\calM \cong \calL$;
\item the set of homomorphisms $\ker(\lambda) \to G(\calL)$ lifting the inclusion $\ker(\lambda) \\\hookrightarrow A$.
\end{itemize}
Under this correspondence $K(\calM) = \ker(\lambda)^\perp/\ker(\lambda)$. In particular, $\calM$ defines a principal polarization on $B$ if and only if $\ker(\lambda)$ is a maximal isotropic subgroup.
\end{prp}

Assume that $\calL$ defines a principal polarization on $A$.  Then $K(\calL^n) = A[n] = \ker ([n]_A)$, where $[n]_A$ is the multiplication map $x\mapsto nx$ in $A$. Applying the previous proposition to $\calL^n$, we obtain

\begin{cor} Assume $(n, \cha~\BK) = 1$. Let $F$ be a maximal isotropic subgroup of $A[n]$ and $\lambda:A\to B= A/F$ be the quotient map. Then  
$B$ admits a principal polarization $\calM$ such that $\lambda^*\calM \cong \calL^n$.
\end{cor}

\subsection{Kummer varieties.} Let $A$ be a principally polarized abelian variety.   Since $h^0(\calL) = 1$, there exists an effective divisor $\Theta$ such that $\calL\cong \calO_A(\Theta)$. Such a divisor $\Theta$ is called a \emph{theta divisor} associated to the polarization. It is defined only up to a translation. One can always choose a  theta divisor satisfying $[-1]_A^*\Theta = \Theta$, which us called a \emph{symmetric theta divisor}. Two symmetric theta divisors differ by a translation $t_a, a\in A[2]$. 

The proof of the following result  over $\BK = \BC$ can be found in \cite{LB}, Chapter IV, \S 8 and in \cite{Ducrohet} in the general case.

\begin{prp} Let  $A$ be a principally polarized abelian variety and $\Theta$ be  a symmetric theta divisor. Then the map $\phi_{2\Theta}:A\to \left| 2\Theta\right|^*$ factors through the projection $\phi:A\to A/\la [-1]_A\ra$  and a morphism  $j: A/\la [-1]_A\ra\hookrightarrow \left| 2\Theta\right|  \cong \BP^{2^g-1}$. If  $A$ is not the product of principally polarized varieties of smaller dimension  and $\cha~\BK \ne 2$, then $j$ is a closed embedding. 
\end{prp}

We assume that $\cha~\BK \ne 2$. The quotient variety $A/\la [-1]_A\ra$ is denoted by $\Kum(A)$ and is called the \emph{Kummer variety} of $A$. In the  projective embedding  $\Kum(A)\hookrightarrow \BP^{2^g-1}$ its degree is equal to $2^{g-1}g!$.
The image of any  $e\in A[2]$ in $\Kum(A)$ is a singular point $P_e$, locally (formally) isomorphic to the affine cone over the second Veronese variety of $\BP^{g-1}$. For any $e\in A[2]$, the image of $\Theta_a := t_e^*(\Theta)$ in $\Kum(A) \subset \BP^{2^g-1}$ is  a subvariety $T_e$ cut out by a hyperplane $2\Theta_e$ with multiplicity $2$. Such a $T_e$ is called a \emph{trope}. Since each $\Theta_e$ is symmetric, the corresponding trope $T_e$ is isomorphic to the quotient $\Theta_e/\la [-1]_A\ra$. 

The configuration of the singular points $P_e$ and the tropes $T_e$ form an abstract symmetric configuration 
$(2^{2g},2^{g-1}(2^g-1))$. This means that each trope contains $2^{g-1}(2^g-1)$ singular points and each singular point is contained in the same number of tropes.

The Kummer variety $\Kum(A)$ admits a resolution of singularities 
\[
\pi:\calK(A)\to \Kum(A)
\]
 with the exceptional locus equal to the union of $E_e = \pi^{-1}(P_e), e\in A[2]$. Each $E_e$ is isomorphic to $\BP^{g-1}$ and the self-intersection $E_e^{g}$ is equal to the degree of the Veronese variety $\nu_2(\BP^{g-1})$ taken with the sign  $(-1)^{g-1}$, that is, the number $(-2)^{g-1}$.

Let $p > 2$ be a prime number and $K$ be a maximal isotropic subgroup in $A[p]$. If 
$p\ne \cha~ \BK$, then $A[p]\cong \BF_p^{2g}$ and the number of such $K$'s is equal to $\prod_{i=1}^g(p^i+1)$. If 
$p = \cha~\BK$, we assume that $A$ is an ordinary abelian variety, i.e. $A[p]_{\text{red}} \cong \BF_p^g$. In this case $F = A[p]_{\text{red}}$ is unique.

\begin{prp}\label{prop1} Let $\lambda:A\to B = A/F$ be the quotient isogeny defined by $\calL^p$.  There exists a symmetric theta divisor $\Theta$ on $A$ and a symmetric theta divisor $\Theta'$ on $B$ such that $\lambda^*\Theta' \in \left| p\Theta\right|$ and $\lambda(\Theta)\in \left| p^{g-1}\Theta'\right| $. Let $D$ be the proper transform of $\lambda(\Theta)$ in $\calK(B)$. Let $m_e$ be  the multiplicity of $\Theta$ at a 2-torsion point $e$. Let
$H$ be the divisor class of the pre-image of a  hyperplane in $\BP^{2^g-1}$ under the composition map $\calK(B)\overset{\sigma}{\to} \Kum(B)\overset{j}{\to}  \BP^{2^g-1}$, where the first map is a resolution of singularities and the second map  is  induced by the map $\phi_{2\Theta}$.
Then
\[
 2D\in \left|  p^{g-1}H-\sum_{e\in A[2]}m_eE_e\right|,
\]
\end{prp}

\begin{proof} As we observed earlier there exists an ample invertible sheaf $\calM$ on $B$ defining a principal polarization such that $\lambda^*\calM \cong \calL^p$. Let $\Theta'$ be a theta divisor on $B$ defined by $\calM$. We have $\lambda^*\Theta'\in \left| p\Theta\right| $ and 
\[
\lambda^*(\lambda(\Theta)) = \sum_{e\in K} t_e^*(\Theta) \equiv p^g\Theta.
\]
Since  the canonical map of the Neron-Severi groups
$\lambda^*:\NS(B) \to \NS(A)$ is injective, we obtain that $p^{g-1}\Theta'$ and $\lambda(\Theta)$ are algebraically equivalent divisors on $B$. Since they are both symmetric divisors, they differ by a translation with respect to a $2$-torsion point $e$. Replacing $\Theta'$ by  $t_a^*(\Theta)$ we obtain the linear equivalence of the divisors. 

It remains to prove the second assertion. Let $\sigma:B' \to B$ be the blow-up of all the $2$-torsion points on $B$. We have a commutative diagram
\[
\xymatrix{
 B'\ar[d]_\sigma\ar[r]^{\phi'}&\calK(B)\ar[d]_{\pi}\\
B\ar[r]^\phi& \Kum(B).
}
\] 
It is clear that $\lambda$ defines a bijection $A[2]\to B[2]$. Since $\lambda$ is a local isomorphism, for any $e\in A[2]$, the multiplicity $m_e$ of $\Theta$ at $e$ is equal to the multiplicity of $\lambda(\Theta)$ at $\bar{e} = \lambda(e)$. Thus 
$2\lambda(\Theta)$ belongs to the linear system $\left| p^{g-1}(2\Theta')-2\sum_{e\in A[2]}m_e \bar{e}\right| $ of divisors in $\left| 2p^{g-1}\Theta'\right|$ passing through the $2$-torsion points $\bar{e}$ with multiplicities $2m_e$. Let $D'$ be the proper transform of $\lambda(\Theta)$ in $B'$. Then
$D' \sim  \sigma^*(p^{g-1}\Theta'-\sum m_e\sigma^{-1}(\bar{e}))$. On the other hand, since $\phi'$ ramifies over each $E_{\bar{e}}$ with multiplicity 2, we have 
\[
 2D'\sim \sigma^*(p^{g-1}2\Theta'-2\sum m_e\sigma^{-1}(\bar{e}))\sim \phi'^*(p^{g-1}\pi^*(H)-\sum_em_eE_{\bar{e}}).
\]
This shows that the proper transform of the image of $2\lambda(\Theta)$ in $\calK(B)$ is linearly equivalent to $p^{g-1}\pi^*(H)-\sum_em_eE_{\bar{e}}$. \end{proof}

\begin{rmr} It is known that a theta divisor on a general principally polarized abelian variety has no singular points at $2$-torsion points. Thus $m_e = 1$ for $2^{g-1}(2^g-1)$ points and $m_e = 0$ at the remaining 2-torsion points. Finally note that $\Pic(\calK(B))$ has no 2-torsion; hence the divisor class $p^{g-1}\pi^*(H)-\sum_em_eE_{\bar{e}}$ has a unique half.
\end{rmr}

\subsection{Theta level structure}\label{heisenberg} The main reference here is \cite{Mu1} (see 
also \cite{Bolognesi}, \cite{LB} where only the case $\BK = \BC$ is considered). Let $A$ be an ordinary abelian variety of dimension $g$ and $\calL \cong \calO_A(\Theta)$ be an ample invertible sheaf defining a symmetric principal polarization. The theta divisor $\Theta$ defines a function 
\[
 q_\Theta:A[2]\to \mu_2,\qquad x\mapsto (-1)^{\textup{mult}_x(\Theta)+\textup{mult}_0(\Theta)}.
\]
This function is a quadratic form whose associated bilinear form is the Weil pairing. We call $\Theta$ {\em even} (resp. {\em odd}) if the quadratic form is even (resp. odd). Recall that the latter means that 
$\# q^{-1}(0) = 2^{g-1}(2^g+1)$ (resp. $\# q^{-1}(1) = 2^{g-1}(2^g-1)$). One can show that $\Theta$ is even if and only if $\textup{mult}_0(\Theta)$ is even. Also, if we normalize the isomorphism $\calL \to [-1]_A^*\calL$ to assume that it is equal to the identity on the fibres over the zero point, then $\Theta$ is even if and only if $[-1]_A^*$ acts as the identity on $\Gamma(\calL)$.

 Let $G(\calL^n)$ be the theta group of $\calL^n$.  A {\em level $n$ theta structure} on $A$ is a choice of an isomorphism of group schemes
\[
\theta:G(\calL^n) \to \calH_g(n),
\]
where $\calH_g(n)$ is the \emph{Heisenberg group scheme} defined by the exact sequence
\[
1\to \BG_m \to \calH_g(n) \to (\BZ/n\BZ)^g\oplus \mu_n^g \to 1.
\]
 It is required that the restriction of $\theta$ to the center is  the identity. 
If $(\lambda,a,b)\in \BK^*\times (\BZ/n\BZ)^g\times \mu_n^g$ represents a point of $ \calH_g(n)$, then the law of composition is 
\[
(\lambda,a,b)\cdot (\lambda',a',b') = (\lambda\lambda'b'(a),a+a',bb'),
\]
where we identify $ \mu_n^g$ with $\Hom((\BZ/n\BZ)^g,\BK^*)$.

A theta level  $n$ structure defines an $n$-level structure on $A$, i.e. an isomorphism of symplectic group schemes
\[
\bar{\theta}: (A[n],e^{\calL^n})\to \bigl((\BZ/n\BZ)^g\times  \mu_n^g,E\bigr),
\]
where
\[
E:(\BZ/n\BZ)^{2g}\times  \mu_n^g\to \BG_m
\]
is the standard symplectic form $((a,b),(a',b')) \mapsto b'(a)/a'(b).$
In particular, $\bar{\theta}^{-1}((\BZ/n\BZ)^{g})$ is a maximal isotropic subgroup in $A[n]$.

The choice of a theta structure of level $n$ defines a representation of the Heisenberg group $\calH_g(n)$ on the linear space $V_n(g) = \Gamma (A,\calL^n)$, called the \emph{Schr\"odinger representation}. In this representation the space $V_n(g)$ admits a basis $\eta_\sigma, \sigma\in (\BZ/n\BZ)^{g},$ such that 
$(\lambda,a,b)\in \calH_g(n)$ acts by sending
$\eta_\sigma$ to $\lambda b(\sigma+a)\eta_{\sigma+a}.$ We will explain 
how to build such a basis from theta functions when we discuss the
 case $\BK = \BC$.

If $n \ge  3$, the map $\phi:A\to \BP(V_n(g)^*)$ given by the complete linear system $\left|\calL^n\right|$ is a closed embedding and the Schr\"odinger representation defines a projective linear representation of the abelian group scheme $(\BZ/n\BZ)^g\oplus \mu_n^g$ in $\BP(V_n(g)^*) \cong \BP^{g^n-1}$ such that the image of $A$ is invariant, and the action on the image is the translation by $n$-torsion points.

Let $\calL$ be a symmetric principal polarization.  The automorphism $[-1]_A$ of $A[n]$  can be extended in a canonical way to an  automorphism $\delta_{-1}$ of $G(\calL^n)$. A theta structure is called \emph{symmetric} if, under the isomorphism $G(\calL^n)\to \calH_g(n)$, the automorphism $\delta_{-1}$ corresponds to the automorphism $D_{-1}$ of $\calH_g(n)$ defined by $(t,a,b)\mapsto (t,-a,b^{-1})$. This defines an action of $D_{-1}$ in $V_n(g)$. 

\smallskip
{\bf From now on we assume that $n > 1$ is odd.}

\smallskip
Since $D_{-1}$ is of order $2$, the vector space $V_n(g)$ decomposes  into the direct sum of two eigensubspaces $V_n(g)^+$ and $V_n(g)^-$ with eigenvalues $1$ and $-1$, respectively.  If $\calL$ is defined by an even theta divisor $\Theta$, then $D_{-1}(\eta_\sigma) = \eta_{-\sigma}$ and we can choose a basis $y_\sigma = \eta_{\sigma}+\eta_{-\sigma}, \sigma\in (\BZ/n\BZ)^g$ in $V_n^+$ and a basis  $z_\sigma =\eta_{\sigma} -\eta_{-\sigma}, \sigma\in (\BZ/n\BZ)^g$ in $V_n^-$\footnote{We order the pairs $\sigma,-\sigma$  to fix the signs of the $z_\sigma$'s.}. In particular,
\[
\dim V_n(g)^{\pm} = (n^g\pm 1)/2.
\]
If $\calL$ is defined by an odd theta divisor $\Theta$, then $D_{-1}(\eta_\sigma) = -\eta_{-\sigma}$ and we have
\[
\dim V_n(g)^{\pm} = (n^g\mp 1)/2 \text.
\]

The two projectivized subspaces form the fixed loci of the projective involution $D_{-1}$. We will call the subspace of dimension $(n^g-1)/2$ the \emph{Burkhardt space} and denote it by $\BP_{Bu}$. The other subspace of dimension $(n^g-3)/2$ we call the \emph{Maschke subspace} and denote it by $\BP_{Ma}$.

 Two different theta structures of level $n$ differ by an automorphism of $\calH_g(n)$ which is the identity on $\BK^*$.  Let $A(\calH_g(n))$ be the group of such automorphisms.  
Let $(\BZ/n\BZ)^{2g}\rtimes \Sp(2g, \BZ/n\BZ)$ be the semi-direct product defined by the natural action of  $\Sp(2g, \BZ/n\BZ)$ on $(\BZ/n\BZ)^{2g}$. There is a natural isomorphism 
\[
(\BZ/n\BZ)^{2g}\rtimes \Sp(2g, \BZ/n\BZ) \to A(\calH_g(n))
\]
defined by sending $(e,\sigma)$ to $(t,u) \mapsto (t[e,u],\sigma(u)).$ 
The group $A(\calH_g(n))$ acts simply transitively on the set of theta structures of level $n$ with fixed even symmetric theta divisor. However, if $n$ is odd, the subgroup of $A(\calH_g(n))$ preserving the set of symmetric structures consists of elements $(0,\sigma)$, hence isomorphic to $\Sp(2g, \BZ/n\BZ) $. This shows that a symmetric theta structure (with fixed $\Theta$) is determined uniquely by the level $n$ structure.

Since the Schr\"odinger representation of $\calH_g(n)$ is known to be irreducible, by Schur's Lemma, the group $A(\calH_g(n))$ has a projective representation in $V_n(g)$. Under this representation, the normal subgroup  $(\BZ/n\BZ)^{2g} \cong \calH_g(n)/\BK^*$ acts via the projectivized Schr\"odinger representation. We will identify  $\Sp(2g, \BZ/n\BZ)$ with the subgroup of $A(\calH_g(n))$ equal to the centralizer of $D_{-1}$. The Burkhardt and the Maschke subspaces  are invariant with respect to the action of the group $\Sp(2g, \BZ/n\BZ)$ in $\BP(V_n(g))$ with the kernel equal to $\la D_{-1}\ra$ and hence  define two  projective representations of $\PSp(2g,\BZ/n\BZ)$ of dimensions $(n^g- 1)/2$ and $(n^g- 3)/2$. 

\subsection{The theta map} We assume that $\cha~\BK \ne 2$ and $n$ is invertible in $\BK$. It is known that the fine moduli spaces of principally polarized abelian varieties of dimension $g$, with a symmetric even or or odd theta structure of odd level $n > 2$, exist. We denote these spaces by $\calA_g(n)^{\pm}$, and by $\calX_g(n)^+$ the corresponding universal families. There is a canonical forgetful morphism 
\[
f_{\pm}:\calA_g(n)^{\pm}\to\calA_g(n)
\]
 to the moduli space of principally polarized abelian varieties of dimension $g$ with level $n$ structure. The fibres are bijective to the set of even (resp. odd) theta divisors, hence the degree of the forgetful map is equal to $2^{g-1}(2^g\pm 1)$.

A theta structure defines a basis  in $V_n(g) = \Gamma(A,\calL^n)$ which is independent of $A$. This defines two $(\BZ/n\BZ\times \mu_n)^g\rtimes \Sp(2g,\BZ/n\BZ)$-equivariant morphisms
\beq\label{exthetamap}
\widetilde{\Th}^{\pm}:\calX_g(n)^{\pm}\to \BP^{n^g-1} = \BP(V_n(g)^*),
\eeq
where the group $(\BZ/n\BZ\times \mu_n)^g$ acts by translations on the image of each $A$.
By composing \eqref{exthetamap} with the zero section $\calA_g(n)^{\pm}\to \calX_g(n)^{\pm}$ we get two $\Sp(2g, \BZ/n\BZ)$-equivariant morphisms, the \emph{even theta map} and the \emph{odd theta map}
\beq\label{thetamap}
\Th^+:\calA_g(n)^{+}\to \BP_{Bu},\quad  \Th^-:\calA_g(n)^{-}\to \BP_{Ma}.
\eeq
Here we use the fact that the value at the origin  of any section of the subspace $V_n(g)^-$  is equal to zero.

\smallskip
 Recall that over $\BC$ the coarse moduli space $\calA_g$ of principally polarized abelian varieties is isomorphic to the orbit space $\calZ_g/\Gamma_g$, where $\calZ_g$ is the Siegel moduli space of complex symmetric $g\times g$-matrices $\tau = X+iY$ such that $Y > 0$, and $\Gamma_g = \Sp(2g,\BZ)$ acts  on $\calZ_g$ in a well-known manner. The moduli space $\calA_g(n)$ is isomorphic to $\calZ_g/\Gamma_g(n)$, where $\Gamma_g(n) = \{M\in \Gamma_g: M-I_{2g} \equiv 0 \mod n\}$. 
It is known that the index of $\Gamma_g(n)$ in $\Gamma_g$ is equal to
$n^{g^2}\prod_{i=1}^g(n^{2i}-1)$, the order of the finite symplectic group $\Sp(2g,\BZ/n\BZ)$. We have a canonical exact sequence
\[
1\to \Gamma_g(2n) \to \Gamma_g(n)\to \Sp(2g,\BF_2) \to 1
\]
defined by the natural inclusion of the groups $\Gamma_g(2n) \subset \Gamma_g(n)$. Comparing the indices with the order of  $\Sp(2g,\BF_2)$, we see that the last map is surjective. It is well-known that the group $ \Sp(2g,\BF_2)$ contains the even and the odd orthogonal subgroups $\Or(2g,\BF_2)^{\pm}$ of indices $2^{g-1}(2^g\pm 1)$. Let $\Gamma_{g}(n)^{\pm}$ be the pre-image in $\Gamma_g(n)$ of the subgroup $\Or(2g,\BF_2)^{\pm}$. Then 
\[
\calZ_g/\Gamma_{g}(n)^{\pm} \cong  \calA_g(n)^{\pm}.
\]
A choice of a symmetric theta structure is defined by a line bundle $\calL$ whose space of sections is generated by a Riemann theta function $\thetafun{m}{m'}(z;\tau)$ with theta characteristic $(m,m')\in (\BZ/n\BZ)^g\times (\BZ/n\BZ)^g$. The even (resp. odd) structure corresponds to the case when $m\cdot m'\equiv 0 \mod 2$ (resp. $\equiv 1 \mod 2$).  A basis of the space $\Gamma(\calL^n)$ is given by the functions
$\thetafun{m+\sigma}{m'}(nz,n\tau)$, where $\sigma\in  (\BZ/n\BZ)^g$. It follows from the standard properties of  the Riemann theta function that 
\beq\label{minus}
\thetafun{m+\sigma}{m'}(-z;\tau) = (-1)^{m\cdot m'}\thetafun{m-\sigma}{m'}(z;\tau).
\eeq 

The theta map \eqref{thetamap} is defined by the theta constants $x_\sigma =\thetafun{m+\sigma}{m'}(0,n\tau)$. They span the space of modular forms of weight $1/2$ with respect to the group $\Gamma(n)^{\pm}$ and some character $\chi: \Gamma(n)^{\pm}\to \BC^*$. It follows from \eqref{minus} that the functions  $y_\sigma = x_\sigma+x_{-\sigma}$ 
(resp.$z_\sigma = x_\sigma-x_{-\sigma}$) are identically zero if $(m,m')$ is odd (resp. even). This shows that the theta maps have the same target spaces as in \eqref{thetamap}.

\begin{prp}\label{SM} Assume $\BK = \BC$. The even theta map
\[
\Th^+:\calA_g(n)^{+}\to \BP_{Bu}
\]
is an embedding for $n = 3$.
\end{prp}

The proof can be found in \cite{SM}, p. 235.

\section{Abelian surfaces}
\subsection{Kummer surfaces.} Now we specialize to the case when $A$ is a principally polarized abelian surface. It is known that $A$ is not the product of two elliptic curves if and only if $\Theta$ is an irreducible divisor.  In this case 
$\Theta$ is a smooth curve   of genus $2$ and $A$ is isomorphic to its jacobian variety $\Jac(\Theta)$. By the adjunction formula, $K_\Theta \cong \calO_{\Theta}(\Theta)$, and so the map $\phi_{2\Theta}$ restricts to the bicanonical map of $\Theta$ onto the corresponding trope of $\Kum(A)$. Let $C$ be a genus $2$ curve and $\Jac^1(C)$ be its Picard scheme of degree $1$. Fix a Weierstrass point $w_0$ to identify $\Jac^1(C)$ with $\Jac(C)$. Then one can take for $\Theta$ the translate of the divisor $W$ of effective divisors of degree $1$, naturally identified with $C$. Under this identification $\Theta$ contains the six $2$-torsion points $w_i-w_0$, where $w_0 = w_1,w_2,\ldots, w_6$ are the six Weierstrass points on $C$. None of them is a singular point of $\Theta$.

Assume $A = \Jac(C)$. In this case $\Kum(A)$ is isomorphic to a quartic surface in $\BP^3$. It has $16$ nodes as singularities and its tropes are conics passing through 6 nodes. The surface $\calK(A)$ is a K3 surface with $16$ disjoint smooth rational curves $E_e, e\in A[2]$. The proper transform of a trope $T$  is a smooth rational curve $\overline{T}$ in the divisor class $\half(H-\sum_{e\in T} E_e)$. 

Assume $A$ is the product of two elliptic curves $F\times F'$. In this case $\Kum(A)$ is the double cover of  a  nonsingular quadric $Q\cong \BP^1\times \BP^1$ branched over the union  $B$ of eight lines, four in each family. Their preimages on $\Kum(A)$ form the union of two sets of four disjoint  smooth rational  curves. The  tropes $T$ on $\Kum(A)$ are the unions of a curve $T_1$ from one set and a curve $T_2$ from another set. Each component of a trope has four $2$-torsion points, one point is common to both components. The proper transform of a trope on $\calK(A)$ is the disjoint union of two smooth rational curves from  the divisor class
$\half(H-\sum_{e\in T_1+T_2}E_e-2E_{T_1\cap T_2})$.

\subsection{Main result I}
We employ the notations of Proposition \ref{prop1}.

\begin{prp}\label{prop2} Assume $A = \Jac(C)$.  Then $\lambda(\Theta)$ is an irreducible curve of arithmetic genus $p^2+1$ with $p^2-1$ ordinary double points.   Its image $D$ in $\Kum(B)$ is a rational curve of arithmetic genus $\half(p^2-1)$ with $\half(p^2-1)$ ordinary double points.
\end{prp}

\begin{proof} We know that $\lambda(\Theta)\in \left| p\Theta'\right|$. Thus $\lambda(\Theta)^2 = 2p^2$ and the first assertion follows from the adjunction formula. Since the isogeny $A\to B$ is a local isomorphism in \'etale topology, the curve $\lambda(\Theta)$ has only ordinary multiple points corresponding to the intersection of an orbit of $K$ with $\Theta$.  Let $\Gamma_a \subset A\times A$ be the graph of the translation map $t_a, a\in K$. It is algebraically equivalent to the diagonal $\Delta_A$ of $A\times A$. Let $C\times C \subset A\times A$ be embedded via the Cartesian product of the Abel-Jacobi map. A point in the intersection $(C\times C)\cap \Gamma_a$ is a pair  of points $(x,y)$ on $C$ such that $[x-y] = a$. By the intersection theory, 
\[
(C\times C)\cdot \Gamma_a  = (C\times C)\cdot \Delta_A = \Delta_C^2 = 2.
\]
Thus, for any nonzero $a\in K$, there exist two  ordered pairs of points on $C$ such that the difference is linearly equivalent to $a$. Since $x-y \sim x'-y'$ implies that $x+y'\sim x'+y \sim K_C$, the two pairs differ by the hyperelliptic involution. If we forget about the order  we get $p^2$ unordered pairs of points in a coset of $K$. This shows that $\lambda(\Theta)$ has $p^2-1$ ordinary double points.

The last assertion follows from Proposition \ref{prop1}.
\end{proof}

\begin{rmr}\label{humbert} 
Note that one can rewrite the formula from Proposition \ref{prop1} in the form
\beq\label{neweq}
\frac{1}{2}(pH-E_1-\ldots-E_6)  = \frac{p+1}{2}H-\frac{1}{2}(H-E_1-\ldots-E_6)-(E_1+\ldots+E_6).
\eeq
Assume that $\Theta'$ is irreducible, i.e. $B = \Jac(C')$ for some curve $C'$. In this case we use $|H|$ to realize $\Kum(B)$ as a quartic surface in $\BP^3$ with 16 nodes. Formula \eqref{neweq} shows that the image $D$ of $\lambda(\Theta)$ on $\Kum(B)$   is cut out  by a surface  $S$ of degree $(p+1)/2$ containing the trope $T$. Since $D$ is a rational curve of arithmetic genus $\half(p^2-1)$ with $\half(p^2-1)$ ordinary double points, the surface $S$ is  tangent to the Kummer surface at $\half(p^2-1)$ points. Humbert proves this fact in \cite{Hum}. In the case $p = 3$ he considers the projection of $D$ to the projective plane from a node of $\Kum(B)$ not lying on the trope $T$. The image of the projection is a rational curve $Q$ of degree 6 that passes through the six nodes of the branch curve, the projections of the double points on the trope $T$. The curve $Q$ is also tangent to the branch curve $B$ at any other intersection point. The six nodes of the branch curve correspond to $6$ parameters of the rational parametrization $\BP^1\to Q$. They are projectively equivalent to the six Weierstrass points defining the original curve $C$. Since the Weierstrass points of the curve $C'$ with $\Jac(C') = B$ are defined by the intersection of $B$ with the osculating conic, this gives a ``solution'' of our problem for $p = 3$. We put it in the quotation marks since it seems impossible to find $Q$ explicitly (we tried!). Note that $Q$ contains $10$ double points, four of them are the projections  of the double points of $D$, the remaining $6$ points are resolved under the double cover. The pre-image of $Q$ under the double cover is equal to the union $D+D'$, where 
\[
D'\in |\frac{1}{2}(9H-12E-E_1-\ldots-E_6)|,
\]
and $E$ is the class of the exceptional curve corresponding to the node from which we project the Kummer surface.

The curves $D,D'$ intersect at $24$ points, $12$ points correspond to the intersection of $Q$ with $B$ outside its singular points, and $12$ other points correspond to six double points of $Q$ that are resolved under the double cover. These peculiar properties of the rational curve $Q$ are impossible to fulfill even with a computer's help. Note that given such a curve $Q$ its pre-image on the Kummer surface defines a curve $D$ linearly equivalent to the divisor from \eqref{neweq}. Its pre-image on the abelian surface $B$ defines a genus $2$ curve whose jacobian variety is isomorphic to $A$. Thus there are $40$ curves $Q$ with the above properties, each defines an isogeny $A\to B$. The problem is that we cannot construct any of them.
\end{rmr}

\subsection{The reducible case}
Note that the principally polarized abelian surface $B = A/F$ could be reducible, i.e. isomorphic to the product $C_1\times C_2$ of two elliptic curves. In this case $B$ is  isomorphic to $\Pic^0(C')$, where $C'$ is a stable curve of genus $2$ isomorphic to the union of two elliptic curves $C_1,C_2$ intersecting transversally at one point.  The theta divisor $\Theta'$ on $B$  is equal to $ (\{0\}\times C_2)\cup (C_1\times \{0\})$, or to its translate by a 2-torsion point. As we mentioned before, the Kummer surface $\Kum(B)$ is a double cover of a nonsingular quadric $Q$. The covering involution $\sigma$  leaves the divisors $H, E_i$ invariant.  However $\sigma$ does not act identically on  the linear system 
$|\frac{1}{2}(pH-E_1-\ldots-E_6)|$. Using the decomposition \eqref{neweq}, we see that the  pre-image of the residual part of a divisor from $|\calO_Q(\frac{p+1}{2})|$ containing the image of the trope $T$ on $Q$ belongs to the linear system $|D|$. The dimension of this linear system is equal to 
\[
\dim |\calO_Q(\tfrac{p-1}{2})| = \frac{1}{4}(p^2+2p-3) < \dim |D| = \half(p^2-1).
\]
The image $D$ of $\lambda(\Theta)$ in $\Kum(B)$ is an irreducible  member of the linear system  $|\frac{1}{2}(pH-E_1-\ldots-E_6)|$ with $\half(p^2-1)$ nodes. It is easy to see that it cannot come from $|\calO_Q(\frac{p-1}{2})|$. Thus  $D$ corresponds to a divisor from $|\calO_Q(p)|$ that passes through the images of the nodes contained in $D$ and splits under the cover.

Consider the case $p = 3$. Let $L_1,\ldots,L_4$ and $L_1',\ldots,L_4'$ be the components of the branch divisor of $\Kum(B)\to Q$, the first four belong to the same ruling. Let $p_{ij} = L_i\cap L_j'$. Their pre-images in $\Kum(B)$ are the $16$ nodes. Without loss of generality we may assume that the curves $E_1,\ldots,E_6$ correspond to the points $p_{12},p_{13},p_{14}, p_{21},p_{31}, p_{41}$.  A curve of bi-degree $(3,3)$ passing through these points splits if and only if it is tangent to each of the  lines $L_2,L_3,L_4,L_2',L_3',L_4'$ at one point and the tangency points are coplanar. Take any rational $4$-nodal curve $R$ of bi-degree $(3,3)$.  These curves depend on $11= 15-4$ parameters. Let $\pi_i:Q\to \BP^1$ be the two rulings on the quadric. The projection $\pi_i:R\to \BP^1$  ramifies over  four points. Let $f_1,\ldots,f_4$ be the corresponding fibres of $\pi_1$ and $ f_1',\ldots, f_4'$ be the corresponding fibres of $\pi_2$. Let 
$f_i\cap C = 2x_i+y_i, f_i' = 2x_i'+y_i'$. Suppose $y_1,y_2,y_3$ lie on a fibre $f_0'$ of $\pi_2$ and $y_1',y_2',y_3'$ lie on a fibre $f_0$ of $\pi_1$. This imposes 4 conditions. Take 
$L_i = f_i, L_i' = f_i', i = 1,2,3,$ and $L_4= f_0, L_4' = f_0'$. Then the  double cover of $Q$ branched over the curves $L_i,L_i', i = 1,\ldots,4$ defines a Kummer surface of a reducible principally polarized abelian variety $B$. The pre-image of the curve $R$ on $\Kum(B)$ splits if the $6$ points 
$L_i\cap L_j', i, j = 1, 2,3,$ are coplanar. This imposes 3 conditions. Counting parameters we see that a curve $R$ always exists. It defines a $4$-nodal divisor $D$ on $\calK(B)$  from the linear system $|\frac{1}{2}(3H-E_1-\ldots-E_6)|$, where $E_1,\ldots,E_6$ correspond to the pre-images of the six points $p_{12},p_{13},p_{14}, p_{21},p_{31}, p_{41}$ in the notation from above.  The pre-image of $D$ on $B'$ is a genus $2$ curve, the jacobian variety of its normalization $C$ is an abelian surface $A$ isogenous to $B$. The six Weierstrass points of $C$ are projectively equivalent to the six coplanar points $p_{ij}$ on the rational curve $R$. This shows that the reducible case realizes, however we do not know how to construct the curve $R$ effectively. 

\subsection{Main result II}
Let us return to the general case $p > 2$. Let $A = \Jac(C)$ and $F$ be a maximal isotropic subspace in $A[p]$. Consider the restriction of the isogeny $\lambda:A\to B = A/F$ to $\Theta$ and compose it with the map $\phi_{2\Theta'}:B\to \Kum(B) \subset \BP^3$ to obtain a map $f:\Theta \to \BP^3$.  Since $\lambda^*(\Theta') \in \left| p\Theta\right|$, the map $f$ is given by a linear system contained in $\left| 2p\Theta\right|$ restricted to $\Theta$. This is the linear system $\left| 2pK_\Theta\right|$. Since $\Theta$ is invariant with respect to the involution   $[-1]_A$, the image of $f$ is equal to the projection  of a rational normal curve $R_{2p}$  of degree $2p$ in $\BP^{2p}=\BP(H^0(\Theta,2pK_\gT)^{ \iota_A})$  from a subspace $L$ of dimension $2p-4$. Let $v_1,\ldots,v_6$ be the images of the six Weierstrass points of $\Theta$ in $R_{2p}$. The divisor  $2\lambda^*(\Theta')$ belongs to $\left| 2p\Theta\right|$ and defines a hyperplane 
$\calH$ in $\BP^{2p}$ which cuts out $R_{2p}$ at $2p$ points containing the points $v_1,\ldots,v_6$. This is because $\lambda^*(\Theta')$ contains $\Theta\cap A[2]$ which we identified with the Weierstrass points. Our main observation is the following.

\begin{thm}\label{main} Let $(z_i,z_i'), i = 1,\ldots, \half(p^2-1),$ be the images on $R_{2p}$ of the pairs of points on $\Theta$ belonging to the same coset of  $K$ and let $\ell_i = \overline{z_i,z_i'}$ be corresponding secant lines of $R_{2p}$. Then the hyperplane $\calH$ intersects the  secants at $(p^2-1)/2$ points which span  a linear subspace  contained in $L\cong \BP^{2p-4}$. The projection of $R_{2p}$ from $L$ maps the points $v_1,\ldots,v_6$ to a  conic $Q$ in $\BP^3$. If $Q$ is irreducible, the double cover of $Q$ branched over the points $v_1,\ldots,v_6$ is a nonsingular curve $C'$ of genus $2$ such that $\Jac(C')\cong B$. If $Q$ is the union of lines then each component has three of the points $v_i$ and the double covers of each line branched along the three points and the intersection point of the line components define two elliptic curves $E$ and $E'$ such that $B = E\times E'$.
\end{thm}

\begin{proof} Assume first that $B \cong \Jac(C')$ for some nonsingular curve $C'$. By Proposition \ref{prop2}, the image of the Veronese curve $R_{2p}$ in $\BP^3$ is a rational curve with $(p^2-1)/2$ ordinary nodes, the images of the points $z_i,z_i'$. This means that each secant $\ell_i$ intersects the center of the projection $L\cong \BP^{2p-4}$.  Since the divisor $\lambda^*(\Theta')$ is the pre-image of a trope in $\BP^3$, the hyperplane $\calH$ must contain the center of the projection $L$. This implies that $L$ intersects the secants $\ell_i$ at the points where $\calH$ intersects them. The image of $R_{2p}$ in $\BP^3$ lies on the Kummer quartic surface $\Kum(B)$ and intersects the trope $T = \Theta'/\la [-1]_B\ra$ at six nodes. The nodes are the images of the Weierstrass points $w_1,\ldots,w_6$. The conic $T$ and the six nodes determine the isomorphism class of the curve $C'$ such that $\Jac(C')\cong B$.

Next assume that $B$ is the product of elliptic curves $F\times F'$. The argument is the same, only this time the image  of $R_{2p}$ lies on the quadric $Q$, which is the image of $\Kum(B)$ in $\BP^3$. The trope $T=\Theta'/\la [-1]_B\ra$ is mapped to the union of two lines $l_1\cup l_2$ intersecting at a point. Each line contains the images of three nodes of $\Kum(B)$. The image of $R_{2p}$ intersects each line at these three points. Again this reconstructs the isomorphism classes of the elliptic curves $F$ and $F'$.
\end{proof}

\section{The case $p = 3$ and $\BK = \BC$}
\subsection{The Burkhardt quartic and the Coble cubic}\label{Sburk}
We specialize the discussion from subsection \ref{heisenberg} to the case $g = 2$ and $n = 3$. In this case we have the theta maps 
\[
  \begin{aligned}
  \Th^+:&\calA(3)^+ \to \BP_{Bu} \cong \BP^4 \\
  \Th^-:&\calA(3)^- \to \BP_{Ma} \cong \BP^3.
  \end{aligned}
\]
According to Proposition \ref{SM}, the first map is an embedding. The second map is an embedding of the open subset of jacobians \cite{Bolognesi}.

It is also known that  the restriction of the map \eqref{exthetamap}
\[
\widetilde{\Th}^{\pm}:\calX(3)^{\pm} \to \BP(V_3(2)) \cong \BP^8
\]
to any fibre $(A,\Theta,\theta)$ defines a closed embedding 
\[
\phi_{\pm}:A \hookrightarrow \BP^8 = \left| 3\Theta\right|^*.
\]
 This embedding is $\calH_2(3)$-equivariant, where $\calH_2(3)$ acts on $A$ via an isomorphism $\bar{\theta}:A[3] \to \BF_3^4$ compatible with the symplectic structures. It acts  on  $\BP^8$ by means of the projectivized Schr\"odinger representation.

We have the following theorem due to A. Coble (for a modern exposition see, for example, \cite{Hunt}, 5.3.1).

\begin{thm}\label{coble} Assume $\cha~\BK \ne 2,3$. Choose the new coordinates in $\BP^8$ as follows.
\[
 \begin{aligned}
 y_0 = \eta_{00}, \ &2y_1 = \eta_{01}+\eta_{02},\ 2y_2 = \eta_{10}+\eta_{20}, \ 2y_3 = \eta_{11}+\eta_{22}, \ 2y_4 = \eta_{12}+\eta_{21}.\\
\ &2z_1= \eta_{01}-\eta_{02}, \ 2z_2= \eta_{10}-\eta_{20}, \ 2z_3= \eta_{11}-\eta_{22}, \ 2z_4= \eta_{12}-\eta_{21}.
\end{aligned}
\]
Then the image $\phi(A)$ is defined by the equations
\beq\label{feq}
\left(\begin{smallmatrix}y_0^2&2(y_1^2-z_1^2)&2(y_2^2-z_2^2)&2(y_3^2-z_3^2)&2(y_4^2-z_4^2)\\
y_1^2+z_1^2&2y_0y_1&2(y_3y_4-z_3z_4)&2(y_2y_4-z_2z_4)&2(y_2y_3-z_2z_4)\\
y_2^2+z_2^2&2(y_3y_4-z_3z_4)&2y_0y_2&2(y_1y_4+z_1z_4)&2(y_1y_3-z_1z_3)\\
y_3^2+z_3^2&2(y_2y_4+z_2z_4)&2(y_1y_4-z_1z_4)&2y_0y_1&2(y_1y_2+z_1z_2)\\
y_4^2+z_4^2&2(y_2y_3+z_2z_3)&2(y_1y_3+z_1z_3)&2(y_1y_2-z_1z_2)&2y_0y_4\end{smallmatrix}\right)\cdot\left(\begin{smallmatrix}\alpha_0\\
\alpha_1\\
\alpha_2\\
\alpha_3\\
\alpha_4\end{smallmatrix}\right) = 0,
\eeq
{\small \begin{eqnarray}\label{twoeq}
z_1\pi_{01}+z_2\pi_{43}+z_3\pi_{24}+z_4\pi_{32}& = &0,\\ \notag
z_1\pi_{43}+z_2\pi_{02}+z_3\pi_{14}+z_4\pi_{13} &=&0,\\ \notag
z_1\pi_{24}+z_2\pi_{14}+z_3\pi_{03}+z_4\pi_{12} &= &0,\\ \notag
z_1\pi_{32}+z_2\pi_{13}+z_3\pi_{12}+z_4\pi_{04} &= &0, \notag
\end{eqnarray}}
where $\pi_{ij} = \alpha_iy_j-\alpha_jy_i$ and  
the vector of the parameters $(\alpha_0,\ldots,\alpha_4)$ is a point $\alpha$ on the {\em Burkhardt quartic} 
\beq\label{burh}
\calB_4:T_0^4+8T_0(T_1^3+T_2^3+T_3^3+T_4^3)+48T_1T_2T_3T_4 = 0.
\eeq
\end{thm}

The vector $\alpha$ depends only the choice of a 3-level structure on $A$ and its coordinates can be identified with explicit modular forms of weight 2 with respect to $\Gamma_2(3)$ (see \cite{FSM2}, p. 253). As we will review below, the
coordinates $T_i$ may be naturally identified with the coordinates $y_i$ in
$\BP_{Bu}$.
One easily notices that the 9 quadratic forms are the partials of a unique cubic form (surprisingly it was missed by Coble). It defines a cubic hypersurface  $\calC_3$ in $\BP^8$, called by the first author, the {\em Coble cubic}. It has a beautiful moduli interpretation in terms of rank 3 vector bundles on the genus 2 curve $C$ (see \cite{Minh}, \cite{Or}).  Thus, the previous theorem expresses the fact that $\phi(A)$ is the singular locus of the Coble cubic $\calC_3$.

The negation involution $[-1]_A$ acts on $\phi_+(A)$ via the projective transformation $\eta_\sigma\mapsto \eta_{-\sigma}$. In  the new coordinates, it is given by
$y_i \mapsto y_i, \ z_j \mapsto -z_j.$
Its fixed locus in $\BP^8$ is the union of two subspaces 
\[
\BP_{Ma} = \{y_0 =\ldots = y_4 = 0\}, \ \BP_{Bu} = \{z_1=\ldots=z_4 = 0\}.
\]
Intersecting $\phi_-(A)$ with $\BP_{Ma}$ we find $6$ points in $A[2]$ lying on $\Theta$. One of them is the origin of $A$. The remaining $10$ points in $A[2]$ form the intersection of $\phi_+(A)$ with $\BP_{Bu}$.  Let us compute this
intersection: Plugging $y_i = 0$ in the equations in Theorem \ref{coble}, we obtain that the parameters $(\alpha_0,\ldots,\alpha_4)$ satisfy the equations
\beq\label{seceq}
\begin{pmatrix}
0&-2z_1^2&-2z_2^2&-2z_3^2&-2z_4^2\\
z_1^2&0&-2z_3z_4&-2z_2z_4&-2z_2z_4\\
z_2^2&-2z_3z_4&0&2z_1z_4&-2z_1z_3\\
z_3^2&2z_2z_4&-2z_1z_4&0&2z_1z_2\\
z_4^2&2z_2z_3&2z_1z_3&-2z_1z_2&0\end{pmatrix}
\cdot\begin{pmatrix}
\alpha_0\\
\alpha_1\\
\alpha_2\\
\alpha_3\\
\alpha_4\end{pmatrix} = 0,
\eeq
As is well-known the coordinates of a non-trivial  solution of a skew-symmetric matrix of corank $1$ can be taken to be the pfaffians of the principal matrices. This gives a rational map 
\[
c_-:\BP_{Ma}\to  \calB_4,
\]
\beq\label{mb}
\begin{aligned}
\alpha_0& =  6z_1z_2z_3z_4 \\
\alpha_1 & = z_1(z_2^3+z_3^3-z_4^3)\\
\alpha_2 & = -z_2(z_1^3+z_3^3+z_4^3)\\
\alpha_3 & = z_3(-z_1^3-z_2^3+z_4^3)\\
\alpha_4 & = z_4(z_1^3+z_2^3-z_3^3)
\end{aligned}
\eeq

We now go back to compute the intersection $\phi_-(A)\cap\BP_{Bu}$: Plugging $z_i = 0$ in \eqref{feq} we obtain that $\alpha$ satisfies the equations
\beq\label{hes}
\begin{pmatrix}y_0^2&2y_1^2&2y_2^2&2y_3^2&2y_4^2\\
y_1^2&2y_0y_1&2y_3y_4&2y_2y_4&2y_2y_3\\
y_2^2&2y_3y_4&2y_0y_2&2y_1y_4&2y_1y_3\\
y_3^2&2y_2y_4&2y_1y_4&2y_0y_1&2y_1y_2\\
y_4^2&2y_2y_3&2y_1y_3&2y_1y_2&2y_0y_4\end{pmatrix}\cdot \begin{pmatrix}\alpha_0\\
\alpha_1\\
\alpha_2\\
\alpha_3\\
\alpha_4\end{pmatrix} = 0,
\eeq
Recall that the Hessian hypersurface  $\Hess(V(F))\subset \BK^m$ of a hypersurface $V(F)$ is defined by the determinant of the matrix of the second partials of $F$. It is equal to the locus of points $x$ such that the polar quadric $P_{x^{m-2}}(V(F))$ of $V(F)$ is singular. The locus of singular points of the polar quadrics  is  the Steinerian hypersurface $\textup{St}(V(F))$. It coincides with the locus of points $x$ such that the first polar $P_x(V(F))$ is singular. One immediately recognizes that  the matrix of the coefficients in \eqref{hes}, after  multiplying the last four rows by $2$, coincides with the matrix of the second partials of the polynomial defining the Burkhardt quartic \eqref{burh}. Thus $\alpha$ is a point on the Steinerian hypersurface of the Burkhardt quartic.  On the other hand, we know that it lies on the Burkhardt quartic. This makes $\calB_4$ a very exceptional hypersurface: it coincides with its own Steinerian. This beautiful fact was first discovered by A. Coble \cite{Coble2}.

The first polar of $\calB_4$ at a nonsingular point is a cubic hypersurface with 10 nodes at nonsingular points of the Hessian hypersurface. Any such cubic hypersurface is projectively isomorphic to the \emph{Segre cubic primal} $\calS_3$ given by the equations in $\BP^5$ exhibiting the $S_6$-symmetry:
\[
Z_0^3+\cdots+Z_5^3 = Z_0+\cdots+Z_5 = 0.
\]
The map from the nonsingular locus of $\Hess(\calB_4)$ to $\calB_4 = \St(\calB_4)$ which assigns to a point $x$ the singular point $\alpha$ of the polar quadric $P_{x^2}(\calB_4)$ is of degree 10. Its fibres are the sets of singular points of the first polars. We will give its moduli-theoretical interpretation in the next section. 

Let 
\[
c_+:\Hess(\calB_4)^{nsg} \to \calB_4, \ (y_0,\ldots,y_4)\mapsto \alpha
\]
be the map given by the cofactors of any column of the matrix of coefficients in \eqref{hes}. 

\begin{thm} The image of the map $\Th^{+}$ is equal to $\Hess(\calB_4)^{nsg}$ and the composition of this map with the map $c_+$ is equal to the forgetful map $\calA_2(3)^+\overset{10:1}{\to} \calA_2(3)$. The composition of the map $\Th^{-}$ with $c_-$ is the  forgetful map $\calA_2(3)^-\overset{6:1}{\to} \calA_2(3)$.
\end{thm}

The first assertion is proved in \cite{vdG} (see also \cite{FSM2}). The second assertion is proved in \cite{Bolognesi} (see also \cite{FSM1}).

\subsection{The $3$-canonical map of a genus 2 curve} Let $(A,\Theta,\theta)$ be a member of the universal family $\calX_2(3)^{-}$. We assume that the divisor $\Theta$ is irreducible, i.e. $A \cong \Jac(C)$ for some smooth genus 2 curve $C\cong \Theta$.  By the adjunction formula, the restriction of the map $\phi_-:A\to \BP^8$ to $\Theta$ is the $3$-canonical map 
\[
\phi_{3K_C}:C\to \left| 3K_C\right|^*\subset \BP^8.
\]
Here the identification of $\left| 3K_C\right|^*$ with the subspace  of $\BP^8 = \left| 3\Theta\right|^*$ is by means of the canonical exact sequence
\[
0 \to \calO_A(2\Theta)\to \calO_A(3\Theta)\to \calO_\Theta(3\Theta) \to 0.
\]
Denote the subspace $\left| 3K_C\right|^* \cong \BP^4$ by $\BP_\Theta^4$. The hyperelliptic involution $\iota_C$ acts naturally on $\BP_\Theta^4$ and  its fixed locus set consists of the union of a hyperplane 
$H_0$ and an isolated point $x_0$. The dual of $H_0$ is the divisor $W = w_1+\cdots+w_6$, where $w_i$ are the Weierstrass points. It coincides with $\phi_-(A)\cap \BP_{Ma}$ and hence 
\[
H_0 = \BP_{Ma}.
\]
The dual of $x_0$ is the hyperplane spanned by the image of the Veronese map $\left| K_C\right|\to \left| 3K_C\right|$. The projection map $C\to H_0$ from the point $x_0$ is a degree 2 map onto a rational norm curve $R_3$ of degree 3 in $H_0$. It is ramified at the Weierstrass points. 

Since $\Theta$ is an odd theta divisor, the image of $\Gamma(\calO_A(2\Theta))$ in $\Gamma(\calO_A(3\Theta))$ is contained in $V_3(2)^- \cong \BC^5$. Thus the image of $V_3(2)^-$ in $\Gamma(\calO_\Theta(3\Theta))$ is the one-dimensional subspace corresponding to the point $x_0$. The projectivization of the image of $V_3(2)^+$ in $\Gamma(\calO_\Theta(3\Theta))$ is the subspace $H_0$. 

Observe that 
\[
\left\{x_0\right\}= \BP_\Theta^{4}\cap \BP_{Bu}.
\]
It is known that the subspace $\BP_\Theta^{4}$ is contained in the Coble cubic $\calC_3$ and $\BP_{Bu}\cap \calC_3$ is equal to the polar cubic $P_\alpha(\calB_4)$ (see \cite{Minh}, Proposition 4.3 and section 5.3). A natural guess is that $x_0 = \alpha$. This turns out to be true.

\begin{lma}\label{point} Let $\alpha = c_-(W) \in \calB_4$. Then, considering $\calB_4$ as
a subset of $\BP_{Bu},$ we have
\[
  x_0 = \alpha.
\]
\end{lma}

\begin{proof} For simplicity of the notation let us denote $\phi_-(A)$ by $A$. Let $I_{A}(2)$ be the subspace of  $S^2V_3(2)$ that consists of quadrics containing  $A$. As we know,  it is spanned by the partial derivatives of the Coble cubic $V(F_3)$. Let $I_{\Theta}(2)$ be the space of quadrics in $\BP_\Theta^{4}$ vanishing on $\Theta$.
The polar map $v\mapsto P_v(F_3)$ defines a $\calH_3(2)\rtimes \la D_{-1}\ra$-equivariant  isomorphism $V_3(2)\to I_{A}(2)$. 

Consider the restriction map 
\[
r:I_{A}(2)\to I_{\Theta}(2).
\]
By \cite{Minh}, Proposition 4.7, this map is surjective. By Riemann-Roch, its kernel $L$ is of dimension $5$. We know that $I_{A}(2) = I_{A}(2)^{+}\oplus I_{A}(2)^{-} = \BC^{5}\oplus \BC^{4}$ with the obvious notation. The subspace $I_{A}(2)^{+}$ is spanned by the four quadrics from \eqref{twoeq}. Obviously they vanish on $\BP_{Ma}\subset \BP_\Theta^{4}$. Since they also contain a non-degenerate curve $\Theta$ they vanish on the whole space $\BP_\Theta^{4}$. Thus $L = I_{A}(2)^+\oplus L^-$, where 
$L^- = L\cap I_{A}(2)^{-}$ is of dimension $1$. In other words, there exists a unique point $x\in \BP_{Bu}$ such that 
the polar quadric $P_x(\calC_3)$ vanishes on $\BP_\Theta^{4}$. It remains to prove that $x_0$ and $\alpha$ both play the role of $x$. 

Recall that the important property of the polar is given by the equality
\[
P_x(\calC_3)\cap \calC_3  = \{c\in \ \calC_3:x\in \BT_c(\calC_3)\},
\]
where $\BT_c(\calC_3)$ denotes the embedded Zariski tangent space.
Since $\BP_\Theta^{4}$ is contained in $\calC_3$, for any $c\in \BP_\Theta^{4}$ we have $\BP_\Theta^{4}\subset \BT_c(\calC_3)$. But $x_0$ belongs to $\BP_\Theta^{4}$, therefore $c\in P_{x_0}(\calC_3)$. This proves that $\BP_\Theta^{4}\subset P_{x_0}(\calC_3)$.

Now consider the polar quadric $P_\alpha(\calC_3)$. It is defined by the quadratic form
\[
(\alpha_0,\ldots,\alpha_4)\cdot M(y,z)  \cdot\left(\begin{smallmatrix}\alpha_0\\
\alpha_1\\
\alpha_2\\
\alpha_3\\
\alpha_4\end{smallmatrix}\right) = 0,
\]
where $M(y,z)$ is the matrix from \eqref{feq}. Restricting the quadric to the subspace $\BP_{Ma}$ we see that it is equal to ${}^t\alpha\cdot M(0,z)\cdot \alpha$, where $M(0,z)$ is the skew-symmetric matrix from \eqref{seceq}. Therefore, it is identically zero on the Maschke subspace, and as above, since it also contains $\Theta$, it must contain the whole $\BP_\Theta^{4}$.
\end{proof}

\subsection{The invariant theory on the Maschke space} Recall that a general point in the Maschke space $\BP_{Ma}$ represents the isomorphism class of a genus $2$ curve together with an odd theta structure. The group $\PSp(4,\BF_3)$ of order $25920$ acts projectively in $\BP_{Ma}$ with the quotient map representing the forgetful map $\calA_2(3)^-\to \calA_2^-$, where $\calA_2^- $ is the moduli space of principally polarized abelian surfaces together with an odd theta divisor. If $C$ is a genus $2$ curve, then the odd theta structure on $\Jac(C)$ is a choice of a Weierstrass point on $C$.

It is known that the projective representation of  $\PSp(4,\BF_3)$ can be lifted to a 5-dimensional linear representation of its central extension $G =\BZ/3\BZ\times \Sp(4,\BF_3)$. The group $G$ act in $\BC^5$ as a complex reflection group with generators of order 3 (Number 32 in Shepherd-Todd's list - see \cite{ST}).  It was first described by Maschke in \cite{Maschke2}. 

We use the coordinates $(z_1,z_2,z_3,z_4)$ in $\BC^4$ corresponding to the projective coordinates in the Maschke space introduced in section 4.1. Let $x\cdot y$ be the standard hermitian dot-product in $\BC^4$.  In these coordinates $G$ is generated by complex reflections
\[
z \mapsto z-\frac{(1-w)z\cdot r}{r\cdot r}r,
\]
where $\eta = e^{2\pi i/3}$ and $r$ is one of the vectors 
\[
(\sqrt{-3},0,0,0),\  (0, \sqrt{-3},0,0),\ (0,\eta^a,\eta^b,\eta^c),\ (\eta^a,0,\eta^b,-\eta^c),
\]
or one of the vectors obtained from these by permutations of the last three coordinates or  multiplied  by $\pm \eta$. 
There are  $40$ different  hyperplanes of fixed points of the complex reflections. Their union is given by the equation
\[
\Phi_{40} = z_1z_2z_3z_4ABCD = 0,
\]
where
\[
\begin{aligned}
A &= (z_2^3+z_3^3+z_4^3)^3-27z_2^3z_3^3z_4^3,\\
B &= (z_1^3-z_2^3+z_3^3)^3+27z_1^3z_2^3z_3^3,\\
C &= (z_1^3+z_2^3-z_3^3)^3+27z_1^3z_2^3z_3^3,\\
D &= (z_1^3+z_3^3-z_4^3)^3+27z_1^3z_3^3z_4^3.
\end{aligned}
\]

\begin{prp}\label{P40}A point $p = (z_1:z_2:z_3:z_4)\in \BP_{Ma}$  is equal to the value  of the theta map  $\Th^-$ at some $(\Jac(C),\Theta)$ with nonsingular $C$ if and only if $\Phi_{40}(z_1,z_2,z_3,z_4) \ne 0.$
\end{prp}

\begin{proof} Consider the rational  map $c_-: \BP_{Ma}\to \BP_{Bu}$ from \eqref{mb}. Its image is the Burkhardt quartic $\calB_4$. The restriction of $c_-$ to any reflection hyperplane defines a rational map $\BP^2\to \BP_{Bu}$ whose image is a plane contained in $\BP_{Bu}$. It is one of the $j$-planes described in \cite{Hunt}. It is shown in loc.cit., Lemma 5.7.3, that 
$\calB_4^0\sm \{40\ j\text{-planes}\}$ is equal to the image of the jacobian locus of $\calA_3(3)$ under the map $c_+:\calA_3(3)^+\to \calB_4$.  One checks that the union of the $j$-planes is equal to the intersection $\calB_4\cap \Hess(\calB_4)$.  Since the Hessian is of degree $10$ and the map $c_-$ is given by quartics, its pre-image in $\BP_{Ma}$ is  a hypersurface of degree $40$. Therefore it coincides with the union of the reflection hyperplanes.

Note that $x\not\in \calB_4\cap \Hess(\calB_4)$ if and only if the polar quadric $P_x^2(\calB_4)$ is singular. The intersection of the cubic polar and the quadric polar is birationally isomorphic to the Kummer surface of a genus $2$ curve $C$ with $\Jac(C)$ corresponding to the point $x$ - see \cite{Coble2}.
\end{proof}

The algebra of $G$-invariant polynomials in $z_1,z_2,z_3,z_4$ was  computed by Maschke \cite{Maschke2}. It is freely generated by polynomials $F_{12}, F_{18}, F_{24}, F_{30}$ of degrees indicated by the subscripts. This shows that 
\[
\BP^3/G \cong \BP(12,18,24,30) \cong \BP(2,3,4,5).
\] 
It follows from Proposition \ref{P40} that the moduli space $\mathcal{M}_2^{\textup{odd}}$ of genus $2$ curves together with a choice of a Weierstrass point is isomorphic to $\BP(2,3,4,5)\sm \{P_{20} = 0\}$, where $P_{20}$ is the polynomial of degree $20$ corresponding to the invariant  $\Phi_{40}^3$.
 This shows that $\mathcal{M}_2^{\textup{odd}}$ is isomorphic  to   $\BP(2,3,4,5)\sm 
V(P_{20})$ for some explicit weighted homogeneous polynomial $P_{20}$ of degree $20$. 

A genus $2$ curve together with a choice of a Weierstrass point can be represented by the equation $y^2+f_5(x) = 0$ for some polynomial of degree 5 without multiple roots. The above discussion suggests that 
the quotient of the open subset of the projectivized space of binary forms of degree $5$ without multiple roots by the affine group $\BC\rtimes \BC^*$  must be isomorphic to $\BP(2,3,4,5)\sm 
V(P_{20})$. This is easy to see directly. Using translations we may choose a representative of an orbit of the form 
\beq\label{f5}
f_5 = x^5+10ax^3+10bx^2+5cx+d.
\eeq
 The group $\BC^*$ acts by weighted scaling
 $(a,b,c,d)\mapsto (t^2a,t^3b,t^4c,t^5d)$. This shows that the orbits of nonsingular quintic forms are parametrized by an open subset of $\BP(2,3,4,5)$. According to G. Salmon \cite{Sa}, p. 230, the discriminant of $f_5$ is equal to 
 \small{$$D = d^4-120abd^3+160ac^2d^2+360b^2cd^2-640bc^3d+256c^5-1440b^3cd^2$$
 $$+2640a^2b^2d^2+4480a^2bc^2d-2560a^2c^4-10080ab^3cd+5760ab^3c^3+3456b^5d$$
 $$+3456a^5d^2-
 2160b^4c^2-11520a^4bcd+6400a^4c^3+5120a^3c^3d-3200a^3b^2c^2.$$}
If we assign to $a,b,c,d$ the weights $2,3,4,5$ respectively, we obtain that $D = P_{20}$ (up to a scalar factor).

\section{An explicit algorithm} 
Let $C$ be a genus $2$ curve and let $F$ be a maximal isotropic subspace in $\Jac(C)[3]$. We would like to find explicitly a stable genus $2$ curve $C'$ such that $\Jac(C')\cong \Jac(C)/F$. 

\subsection{$\BK = \BC$} We start with the complex case. Unfortunately,  we do not know how to input explicitly the pair $(C,F)$. Instead we consider $C$ with an odd theta structure. It follows from Proposition \ref{P40} that the isomorphism class of such a structure  $(\Jac(C),\Theta,\bar{\theta})$ is defined by a point $p$ in $\BP_{Ma}$ not lying on the union of the reflection hyperplanes. The theta structure defines a maximal isotropic subspace $F$ in $\Jac(C)[2]$. Two points $p$ and $q$ define the same maximal isotropic subspace if and only if they lie in the same orbit with respect to some stabilizer subgroup of $\PSp(4,\BF_3)$: the maximal isotropic subspace $F_0$ of the symplectic space $\BF_3^4$ that consists of vectors with the first two coordinates equal to zero. Its  stabilizer subgroup is a maximal subgroup of $\PSp(4,\BF_3)$ of index 40 isomorphic to $(\BZ/3\BZ)^3\rtimes S_4$. 

\smallskip\noindent
\emph{Step 1}: Evaluating the Maschke fundamental invariants $F_{12}, F_{18}, F_{24}, F_{30}$ at $p$ we can find an equation $y^2+f_5(x) = 0$ of $C$, where $f_5$ is as in \eqref{f5}.

\smallskip\noindent
\emph{Step 2}: Next we consider the point $\alpha = c_-(p)\in \calB_4$. From Lemma \ref{point} we know that the span of $\BP_{Ma}$  and $\alpha$ in $\BP^8$ is the space $\BP_\Theta^4$ where the tri-canonical image of $C$ lies.  We know from the proof of Theorem \ref{main} that there are two pairs of points $(x,y), (x',y')$ on $C$ such that $x-y = x'-y' = e$. The group $\BF_3^4$ acts on $\BP^8$ via the Schr\"odinger representation. Take a point $e\in F_0$ that corresponds to some point in $F$ under the theta structure. Intersecting $\BP_\Theta^4$ with its image $\BP_\Theta^4+e$ under  $e$, we get a plane $\Pi_e$. It is clear that $\BP_\Theta^4$ intersects $\BP_\theta^4+e$  along the plane $\Pi_e$ spanned by the secant lines $\overline{x,y}$ and $\overline{x',y'}$. The hyperelliptic involution acting in $\BP_\theta^4$ switches the two lines. In particular,  the lines intersect at a point $p_e\in H_0 = \BP_{Ma}$ (cf. \cite{Bolognesi}, Lemma 5.1.4). Note that replacing $e$ with $-e$ we get the same pair of secants. In this way we obtain 8 secant lines $\overline{x_e,y_e}, \overline{x_e',y_e'}$, each pair corresponds to the pair $(e,-e)$ of 3-torsion points from $F$. Thus the plane  $\Pi_e$  intersects  $C$ at $4$  points $x_e,y_e,x'_e,y'_e$. They define a pair of concurrent secants. 

\smallskip\noindent
\emph{Step 3}: Using $\alpha$ we can write down the equations defining  the abelian surface $A = \Jac(C)$ in $\BP^8$ as given in  \eqref{feq}. Intersecting $A$ with $\BP_\Theta^4$ we find the equations of the tri-canonical model of $C$ in $\BP_\Theta^4$. They are given by the restriction of the quadrics containing $A$ to $\BP_\Theta^4$ (\cite{Minh}, Proposition 4.7). 

\smallskip\noindent
\emph{Step 4}: Next we project $C$ from $\alpha$ to $\BP_{Ma}$. The image is a rational normal curve $R_3$. The image of the 4 pairs of concurrent secants $l_1,\ldots,l_4$ is a set of 4 secants of $R_3$. 

\smallskip\noindent
\emph{Step 5}: Now we need to locate the $6$-tuple of points $\{x_1,\ldots,x_6\}$ on $R_3$ corresponding to the Weierstrass points of $C$. This
is the branch locus of the map $C\to R_3$, which is computed explicitly
by choosing a rational coordinate on $R_3$ (note that we do not need the
coordinate of any of the $x_i$'s by itself).

\smallskip\noindent
\emph{Step 6}: Identifying $R_3$ with $\BP^1$ let us consider the Veronese map $R_3\to R_6 \subset \BP^5$. Let $y_1,\ldots,y_6$ be the images of  the six points $x_1,\ldots,x_6$ and $\ell_1,\ldots,\ell_4$ be the secants defined by the images of the 4 pairs of points defining the secants $l_1,\ldots,l_4$ of $R_3$. The points $y_1,\ldots,y_6$ span a hyperplane $H$ in $\BP^5$. 

\smallskip\noindent
\emph{Step 7}:  This is our final step. Following the proof of Theorem \ref{main}, we intersect the four secants $\ell_1,\ldots,\ell_4$ with $H$. The four intersection points span a plane $\pi$. We project $H$ from $\pi$ to $\BP^2$. The images of the points $y_1,\ldots,y_6$ lie on a conic and determine a stable genus 2 curve $C'$ with $\Jac(C') \cong \Jac(C)/F$.

\subsection{The case $\cha~\BK=3$}\label{Salgo3} Recall (see the discussion before Proposition \ref{prop1}) that
$F = \jac(C)[3]_{\text{red}}\cong (\BZ/3\BZ)^2$. The algorithm we present below
gives the curve $C'$ such that $\jac(C')=\jac(C)/F$. This
construction is
known to be the inverse of the Frobenius map on $\calA_2(\BK)$ (this is
seen by considering the quotient of $\jac(C)$ by the Weil pairing dual
of the group $\jac(C)[3]_{\text{red}}$, which is the group scheme isomorphic to $\mu_{3,\BK}^2$).

 First let us recall  an explicit algorithm for finding $3$-torsion points on $\Jac(C)$ \cite{CF}. Let $w_1,\ldots,w_6$ be the Weierstrass
points of $C$. Fix one of them, say  $w = w_1$, i.e. choose to define $C$ by the equation $y^2t^3-f_5(x,t) = 0$ in $\BP^2$. The  plane quintic model $C_0$ has a triple cusp at $(t,x,y) = (0,0,1)$.  It is the projection of the quintic curve $C$ in $\BP^3$ embedded by the linear system 
$\left| 2K_C +w\right|$ from any point, not on $C$, lying on the ruling of the unique quadric containing $C$ which cuts out the divisor $3w$ on $C$. The pencil of lines through the singular point of $C_0$ cuts out the linear system $\left| K _C\right| +3w$.

A plane cubic with equation $yt^2-f_3(x,t) = 0$ intersects $C_0$ at $6$ nonsingular points 
$p_1,\ldots,p_6$
and at the point  $(0,0,1)$ with multiplicity $9$. This implies that $p_1+\cdots +p_6$ is linearly equivalent to $3K_C$. Using the well-known description of $\Jac(C)$ in terms of the symmetric square of $C$,  we see that $[p_1+p_2]\oplus [p_3+p_4]= -[p_5+p_6]$ in the group law on  $\Jac(C)$, where $[p+q]$  is the divisor class of the divisor $p+q-K_C$. Here  $-[p+q]$ is equal to $[p'+q']$, where $p\mapsto p'$ is the hyperelliptic involution $(t,x,y)\mapsto (t,x,-y)$. 

Let us choose the coordinates so that the affine piece of $C_0$ is given by
$y^2=x^5+\sum_{i=0}^4 b_ix^i$.  Replacing $x$ with $x+\frac{b_4}{5}$ we may assume that $b_4 = 0$. It is known that $\Jac(C)$ is an ordinary abelian variety if and only if the Cartier-Manin matrix 
$$A = \begin{pmatrix}b_2&b_1\\
1&b_4\end{pmatrix}$$
is nonsingular - see \cite{Yui}, Theorem 3.
Since we assumed that $b_4 = 0$, this is equivalent to 
\[
b_1 \ne 0.
\]

We have to find coefficients $a,d_0,d_1,d_2,c_0,c_1$ which
solve the equation
\[
  (x^3+d_2x^2+d_1x+d_0)^2-a(x^5+\sum_{i=0}^3 b_ix^i)-(x^2+c_1x+c_0)^3=0.
\]
Equating coefficients at powers of $x^i$ we find
\begin{eqnarray*}
i = 5: \quad d_2& = &-a,\\ \notag
i= 4: \quad d_1&=&a^2,\\ \notag
i= 3: \quad c_1^3&=&-ab_3+a^3-d_0,\\ \notag
i= 2: \quad d_0 &=&b_2-a^3,\\ \notag
i=1: \quad b_1&=&-ad_0,\\ \notag
i = 0:\quad c_0^3&=&d_0^2-ab_0.
\end{eqnarray*}
The first, the second, and the fourth equations eliminate $d_0,d_1,d_2$. The fifth equation gives a quartic equation for the variable $a$
\beq\label{a}
X^4-b_2X-b_1=0.
\eeq
Note that equation \eqref{a} is separable if and only if  $b_1 \ne  0$.  

It is easy to find the roots of this equation since the splitting field
of the resolvent polynomial is an Artin-Schrier extension of $\BK$.
Finally, the third and the last equations give 
\[
  c_1^3 = a^3-ab_3+\frac{b_1}{a}, \quad 
c_0^3 =  \frac{b_1^2}{a^2}-ab_0.
\]
Each solution $a$ of  \eqref{a} defines a quadratic equation $x^2+c_1x+c_0$ and thus also the pair $\{x_a,x_a'\}$ of its roots. These are the $x$-coordinates (or, equivalently, the orbits with respect to the hyperelliptic involution $\iota$) of a pair of points $p_1,p_2$ such that $p_1+p_2-H$ is a $3$-torsion divisor class in $\Jac(C)$.  Since $\Jac(C)[3]_{\text{red}} \cong (\BZ/3\BZ)^2$, we have four distinct pairs of roots. Consider the pairs of roots as points on the image $R_6$ of $C$ under the map given by the linear system $\left|6K_C\right| ^\iota$. They define the four secants from the proof of Theorem \ref{main}. Now we finish as in Step 7 from the case $\BK = \BC$.

\end{document}